\documentclass[12pt]{article}
\usepackage{setspace}
\usepackage[utf8]{inputenc}\usepackage{authblk}
\usepackage{graphicx}
\usepackage{amsmath}
\usepackage{amssymb,amsthm}
\usepackage{geometry}
\usepackage{cases}
\usepackage[hidelinks]{hyperref}
\usepackage{cleveref}

\graphicspath{ {images/} }
\newtheorem{theorem}{Theorem}

\newtheorem{cor}{Corollary}
\newtheorem{prop}{Proposition}

\theoremstyle{definition}
\newtheorem{rk}{Remark}

\DeclareMathOperator*{\supess}{ess\,sup}
\DeclareMathOperator*{\supp}{supp}

\begin{document}

\def\R{\mathbb{R}}
\def\He{\mathbb{H}}
\def\C{\mathbb{C}}
\def\g{\rm{\bf g}}\def\ric{{\text {Ric}}}
\def\a{\alpha}
\def\e{\epsilon}\def\s{\sigma} 
\def\si{\sigma}
\def\ga{\gamma}
\def\la{\lambda}
\def\k{\kappa}
\def\ds{\displaystyle}
\def\ni{\noindent}
\def\avint{\mathop{\,\rlap{\bf--}\!\!\int}\nolimits}
\def\ba{\begin{aligned}}
\def\ea{\end{aligned}}
\def\ol{\overline}
\def\wt{\widetilde}
\def\b{\beta}\def\ab{{\overline\a}}
\def\avgdj{\avint_{(D_{j+1}\setminus D_{j})\setminus E_\tau}}
 \def\avgdJ{\avint_{(D_{J+1}\setminus D_{J})\setminus E_\tau}}
 \def\avgdJo{\avint_{(D_{J_1+1}\setminus D_{J_1})\setminus E_\tau}}
 \def\dJ{\int_{(D_{J+1}\setminus D_{J})\setminus E_\tau}}
 \def\dj{\int_{(D_{j+1}\setminus D_{j})\setminus E_\tau}}
 \def\wphi{\wt{\phi}_{\e,r}}
 \def\mc{\mathcal}\def\inj{{\text {inj}}}
 \def\wh{\widehat}
 \def\AA{\max\{A_0,A_\infty\}}
\def\be{\begin{equation}}\def\ee{\end{equation}}\def\nt{\notag}\def\bc{\begin{cases}}\def\ec{\end{cases}}\def\ba{\begin{aligned}}
\def\ea{\end{aligned}}
\def\qed{{\setlength{\fboxsep}{0pt}\setlength{\fboxrule}{0.2pt}\fbox{\rule[0pt]{0pt}{1.3ex}\rule[0pt]{1.3ex}{0pt}}}}
\def\QEDopen{{\setlength{\fboxsep}{0pt}\setlength{\fboxrule}{0.2pt}\fbox{\rule[0pt]{0pt}{1.3ex}\rule[0pt]{1.3ex}{0pt}}}}
\def\na{\sum_{k=0}^{2N-1}}
\def\ia{\sum_{k=2N}^{\infty}}
\def\aa{\sum_{k=0}^{\infty}}\def\Tfc{(Tf)^\circ}
\def\udl{\underline}\def\llceil{\left\lceil}\def\rrceil{\right\rceil}
\def\bfe{{\bf e_1}}\newcommand\fH[1]{\sbox0{#1}\dimen0=\ht0 \advance\dimen0 -1ex
  \sbox2{\'{}}\sbox2{\raise\dimen0\box2}%
  {\ooalign{\hidewidth\kern.1em\copy2\kern-.5\wd2\box2\hidewidth\cr\box0\crcr}}}
\def\S{\mathbb{S}}\def\SS{\mathcal{S}}
\centerline{\bf \large{ Moser-Trudinger inequalities: from local to global}}
\def\N{{\mathbb N}}
\font\twelvemi=cmmi12 at 13pt\font\elevenmi=cmmi11 at 8 pt
\renewcommand{\chi}{\raisebox{.13\baselineskip}{\hbox{\twelvemi\char31}}}
\renewcommand{\gamma}{{\raisebox{.08\baselineskip}{\hbox{\twelvemi\char13}}}}
\newcommand{\sgamma}{{\raisebox{.08\baselineskip}{\hbox{\elevenmi\char13}}}}
\def\lan{\big\langle}\def\ran{\big\rangle}
\vskip 0.1in
\centerline{\bf{Luigi Fontana, Carlo Morpurgo,  \ Liuyu Qin\footnote{The third author was supported by the National Natural Science Foundation of China (12201197) \;\;\;\; 
\smallskip MSC Class: 46E36; 26D10}}}
\def\gn{\gamma_n}\def\sgn{\sgamma_{\!n}}
\def\bu{$\bullet\;$} 
\def\nn{{\frac{n}{n-1}}}\def\un{{\frac{1}{n-1}}}

%
%
%
\begin{abstract} 

Given a general complete Riemannian manifold $M$, we introduce the concept of ``local Moser-Trudinger inequality on $W^{1,n}(M)$''. We show how the validity of the Moser-Trudinger inequality can be extended from a local to a global scale under additional assumptions: either by assuming the validity of the Poincar\'e inequality, or by imposing a stronger norm condition. We apply these results to Hadamard manifolds.
The technique is general enough to be applicable also in sub-Riemannian settings, such as the Heisenberg group. 

\end{abstract}
\numberwithin{equation}{section}

 \section{Introduction}

\smallskip

Let $(M,g)$ be a complete, $n-$dimensional  Riemannian manifold, without boundary, and let $\mu$ be  the associated Riemannian measure. We say that a \emph{local Moser-Trudinger inequality holds on $W^{1,n}(M)$} (or ``on $M$'', in short),  if there are constants  $\tau, C_0>0$ such that for any bounded,  measurable $E\subseteq M$ with $\mu(E)\le \tau$
\be\int_{E} e^{\sgn|u|^\nn} d\mu\le C_0,\quad  u\in W^{1,n}(M),\;\|\nabla u\|_n\le 1,\;  \supp u\subseteq E,\label{LMT}\tag{LMT}\ee
where $\gamma_n=n\omega_{n-1}^{\frac1{n-1}}$, and   $\omega_{n-1}$ is the surface measure of the standard $(n-1)$-dimensional,  unit sphere $S^{n-1}$ (for the definition of the usual Sobolev space $W^{1,n}(M)=W_0^{1,n}(M)$  see [H]). The constant $\gamma_{n}$ is the sharp exponential constant in the classical Moser-Trudinger inequality on bounded domains of $\R^n$ [M], or on compact Riemannian manifolds [F]. Needless to say, the prime example of a manifold where \eqref{LMT} holds is $\R^n$ (see Observation 3 below).

\smallskip
We make now a few important observations.

\smallskip
\ni {\bf Observation 1.} 
    On an arbitrary complete Riemannian manifold $M$, given $x_0,\in M$ and $0<r_0<\inj_{x_0}(M)$, the injectivity radius at $x_0$, one can consider the ``Moser family'' of  functions $v_\e\in W^{1,n}_0(B(x_0,r_0))$ defined as 
\be v_\e(y)=\bc \log\dfrac1\e & \text{ if } d(x_0,y)\le \e r_0\cr
\log\dfrac {r_0}{d(x_0,y)} &\text {if } \e r_0<d(x_0,y)\le r_0\cr
0 & \text {if } d(x_0,y)>r_0\ec
\ee
and show, by a standard argument (see e.g. [F]),  that for any $\gamma>\gamma_{n}$, if $u_\e=v_\e/\|\nabla v_\e\|_n$ then
\be \int_{B(x_0,r_0)} e^{\sgamma|u_\e|^\nn} d\mu\ge C \mu(B(x_0,\e r_0))^{1-\sgamma/\sgn}\to+\infty,\qquad \text {as }\e\to0.\ee
 This shows  that $\gamma_n$ in \eqref{LMT} is sharp, i.e. it cannot be replaced by a larger constant.

\medskip 

\ni {\bf Observation 2.}  For \eqref{LMT} to hold it is necessary that the volumes of geodesic balls of fixed radii are not collapsing, i.e. $\mu(B_x(r))\ge C(r)>0$ for all $r>0$ and $x\in M$, where $B_x(r)$ is the geodesic ball centered at $x$ and with radius $r$, and where $C(r)$ is a constant which depends  on $n,\tau,C_0,r$, but is  independent of $x$ (see Prop. 2 in Section 3 for a proof). It is well known that such condition follows if $\inj(M)>0$, where $\inj(M)$ denotes the injectivity radius of $M$. Indeed, Croke [C] proved that if $\inj(M)>0$ then  $\mu(B(x,r))\ge C_n\min\{ (\inj(M)/2)^n, r^n\}$, where $C_n$ is an explicit constant depending on $n$ (see also [CH] for related comments).


\medskip\ni {\bf Observation 3.}
The validity of \eqref{LMT} for any $\tau>0$ is guaranteed provided one can find a measurable function $k(x,y)$ on $M\times M$ such that for all $u\in C_c^\infty(M)$ 
\be |u(x)|\le \int_M k(x,y)|\nabla u(y)|d\mu(y)\label{rep}\ee
where $|Z|^2=g(Z,Z)$, for $Z\in T_yM$, and where $k$ satisfies the following growth conditions for some $C_1,C_2>0$, $t_0>0$,  and $\delta>1$:
\be (k_x)^*(t)\le \gn^{-\frac{n-1}n} t^{-\frac{n-1}n}\big(1+C_1(1+|\log t|)^{-\delta}\big),\qquad  0<t\le t_0, \quad {\hbox{ a.e. }}x\in M \label{kstar1}\ee
\be (k^y)^*(t)\le C_2  t^{-\frac{n-1}n},\qquad 0<t\le t_0, \quad {\hbox{ a.e. }}y\in M.\label{kstar2}\ee
Here $k_x$ and $k^y$ denote the usual $x-$section and $y$-section of $k(x,y)$, and for a measurable function $f$ on $M$ we define its nonincreasing rearrangement as  
\be f^*(t)=\inf\big\{s>0:\; \mu\{x:|f(x)|>s\})\le t\big\},\qquad  t>0\ee
(we will implicitly assume that all level sets of rearranged functions have finite measure).

If the estimates in \eqref{kstar1} and \eqref{kstar2} hold then they also hold for $0<t\le \tau$, any $\tau>0$, with $C_1,C_2$ depending on $\tau$, and    if $T$ is the integral operator with kernel $k(x,y)$, then there is $C$ (depending on $\tau$) such that for any $E,F$ measurable with $\mu(E),\;\mu(F)\le \tau$
\be \int_E e^{\sgn |Tf|^{\nn}}d\mu\le C\big(\mu(E)+\mu(F)\big) ,\qquad f\in L^n(F),\; \|f\|_n\le1.\label{ad}\ee
This result follows from   [FM2, Thm. 1], where a more general statement on arbitrary  measure spaces is given (see also [FM2, Cor. 2], where it is assumed that the estimates in \eqref{kstar1} and \eqref{kstar2} hold for $t>0$).
  Clearly,  \eqref{rep} and \eqref{ad}, with $E=F$,   imply \eqref{LMT}, and in fact a slightly stronger statement, due to the presence of $\mu(E)$ on the RHS of \eqref{LMT}.
 
  The standard way to estimate $u$ in terms of its gradient, as in \eqref{rep}, is via the Green function $G(x,y)$ for the Laplace-Beltrami operator, and the formula
 \be u(x)=-\int_M G(x,y)\Delta u(y)d\mu(y)=\int_M\lan\nabla_y G(x,y),\nabla u(y) \ran d\mu(y),\qquad u\in C_0^\infty(M), \label{rep1}\ee
where $\lan Z,W\ran=g(Z,W)$, for $Z,W\in T_y M$, in which case we can take $k(x,y)=|\nabla_y G(x,y)|$ (which in general is not symmetric). Hence, the problem of  establishing \eqref{LMT} reduces, typically,  to the problem of finding 1)  a sharp \emph{uniform} asymptotic estimate for $|\nabla_y G|$ near the diagonal and 2) a  \emph{uniform} estimate for $|\nabla_y G(x,y)|$ for large distances. Uniformity here means that the constants involved in such estimates, and ultimately $C_1, C_2$ in \eqref{kstar1}, \eqref{kstar2},  are depending only on~$M$. Even assuming $\inj(M)>0$ one typically needs more assumptions on the curvature, in order to obtain such uniform estimates.

The validity of \eqref{LMT} on $\R^n$ is  easily checked, since estimate \eqref{rep} holds with $k(x,y)=\omega_{n-1}^{-1}|x-y|^{1-n}$, obtained from 
\be u(x)=-{\frac1{(n-2)\omega_{n-1}}}\int_{\R^n}|x-y|^{2-n}\Delta u(y)d\mu(y),\ee 
and clearly \eqref{kstar1}, \eqref{kstar2} are satisfied by this kernel.
\def\vol{{\text {vol}}}

\medskip \ni {\bf Observation 4.}  The  statement in \eqref{LMT} holds if $E\subseteq \Omega_0$, where $\Omega_0$ is any \emph{fixed, bounded} open set with smooth 
boundary $\partial \Omega_0\neq \emptyset$.  This statement does not appear in the literature in this exact form, but it is an easy consequence of a  representation formula as in \eqref{rep1}, using the Green function $G_{\Omega_0}^{}$ for the Dirichlet problem on $\Omega_0$, combined with the estimate
\be |\nabla_y G_{\Omega_0}^{}(x,y)|\le \omega_{n-1}^{-1}d(x,y)^{1-n}\big(1+H_{\Omega_0}^{} d(x,y)\big),\qquad x,y\in \Omega_0,\label{nablag}\ee
where $H_{\Omega_0}^{}$ is some constant depending on $\Omega_0$. Such estimate can be derived as in [F], following essentially the argument given by Aubin in [Au1] (see also [Au2, Thm. 4.17]). It's easy to check that \eqref{nablag} implies \eqref{kstar1}, \eqref{kstar2}, where the rearrangements are taken w.r. to $\Omega_0$, and  \emph{with constants $C_1,C_2$ depending on $\Omega_0$}; as a consequence, \eqref{LMT} also holds, \emph{with a constant $C_0$ depending on $\Omega_0$}. It is well known that on any complete  noncompact manifold $M$ (without boundary)  one can find a nested sequence of open bounded sets $\Omega_k$ with smooth boundary, such that $\bigcup_k\Omega_k=M$, and construct the Green function $G(x,y)$  on $M$ by approximation, using the Green functions $G_{\Omega_k}$ (see for ex. [Li]). However, without further assumptions on $M$ there is no guarantee that the Green function so obtained  will satisfy a   uniform bound as in \eqref{nablag} for $d(x,y)\le 1$, or even the uniform bound  $|\nabla_y G(x,y)|\le C$ for $d(x,y)\ge1$.

\medskip \ni {\bf Observation 5.}  If $M$ is compact and $\tau=\mu (M)$, then \eqref{LMT} cannot hold as stated, and it reduces to the classical result in  [F] under the additional assumption that $u$ has mean 0.  We do not know  whether or not,  for $M$ compact,  \eqref{LMT} holds for  $\tau<\mu(M)$, or even for some $\tau$ small enough, without assuming $u$ with 0 mean.

\smallskip

\def\Co{C_\Omega}

\bigskip
We now discuss the validity of the Moser-Trudinger inequality on a global scale. Given an open set $\Omega\subseteq  M$ with infinite measure,  we say that the \emph{classical Moser-Trudinger inequality holds on $W_0^{1,n}(\Omega)$}, if there is $C_\Omega>0$ such that
\be\int_\Omega \exp_{ n-2}\big[\gamma_n|u|^\nn\big] d\mu\le \Co,\qquad  u\in W_0^{1,n}(\Omega),\quad \|\nabla u\|_n\le 1,\label{CMT}\tag{CMT}\ee
where for $m\in\N$
\be \exp_m[t]=e^t-\sum_{k=0}^m\frac{t^k}{k!}.\notag\ee

If $M=\R^n$ and  $\Omega$ has infinite inradius, then a dilation argument shows that \eqref{CMT} fails, even if $\gamma_n$ is replaced by any $\gamma>0$. In [BM] Battaglia and Mancini proved that  \eqref{CMT} on a given open set $\Omega\subseteq \R^n$,  is in fact equivalent to the validity of the Poincar\'e inequality on $W^{1,n}_0(\Omega)$, extending earlier 2-dimensional results by Mancini and Sandeep [MS].

 In this note we establish this equivalence on general complete Riemannian manifolds satisfying \eqref{LMT}.

\smallskip
 Given   $(M,g)$ a  $n$-dimensional complete  Riemannian manifold and an open set $\Omega\subseteq M$, we say that the \emph{Poincar\'e inequality holds  on $W_0^{1,n}(\Omega)$} if there is $C_\Omega$ such that
\be \|u\|_n\le\Co  \|\nabla u\|_n,\qquad u\in W_0^{1,n}(\Omega).\label{Pn}\tag{Pn}\ee

\begin{theorem}\label{1}  Let  $(M,g)$ be an  $n$-dimensional complete  Riemannian manifold such that the local Moser-Trudinger inequality  \eqref{LMT} holds on $M$. Then, for any open $\Omega\subseteq M$, the  Moser-Trudinger inequality \eqref{CMT} holds on $\Omega$ if and only if the Poincar\'e inequality \eqref{Pn} holds on $W_0^{1,n}(\Omega).$ 

\end{theorem}

\medskip

\ni As we mentioned earlier, the obvious example where \eqref{LMT} holds is $\R^n$, 
thus, we instantly recover Battaglia-Mancini's result on $\R^n.$

More general  manifolds where \eqref{LMT} holds,  are Hadamard manifolds with bounded curvature. As a consequence, we can extend Battaglia-Mancini result on such manifolds:

\begin{theorem}\label{2}  Let  $(M,g)$ be an  $n$-dimensional complete, simply connected   Riemannian manifold, and with sectional curvature $K$  satisfying $-b^2\le K\le 0$. Then \eqref{LMT} holds, and  for any open $\Omega\subseteq M$, the  classical Moser-Trudinger inequality \eqref{CMT} holds on $W_0^{1,n}(\Omega)$ if and only if the Poincar\'e inequality \eqref{Pn} holds on $W_0^{1,n}(\Omega).$ 

\end{theorem}
\def\ric{{\text {Ric}}}

We point out here that in [YSK, Theorem 1.1] it is claimed that a local form of the MT inequality (which implies \eqref{LMT})  holds on an arbitrary Hadamard manifolds, i.e. without imposing a lower bound  on $K$. Unfortunately, there appears to be a gap in the proof (see comments in the ``Second proof of Theorem 2'' below). We do not know at this point if \eqref{LMT} holds under the sole condition that $K\le 0$.

Note that under the restriction $K\le -a^2<0$ the Poincar\'e inequality is valid on $W_0^{1,p}(M)$  for all $p\ge 1$ (see [S, Thm. 5.4]), and when $-b^2\le K\le -a^2<0$ and $n\ge 3$,  then \eqref{CMT} holds on $W^{1,n}(M)$, in fact it also  holds a on higher order Bessel potential spaces [BS].

We will present two proofs of Theorem 2. The first uses sharp bounds for the gradient of the Green function, essentially all contained in [BS],  in the case $n\ge3$. The second proof  is based on the method used in [YSK], which works out under the additional assumption $K\ge -b^2$, for any $n\ge2.$

We will not attempt here to characterize, or even discuss, the validity of the Poincar\'e inequality on open sets of general  manifolds satisfying $-b^2\le K\le 0$. We will, however, summarize some of the known results on $\R^n$ in Section 3.
Likewise, in this paper we will not attempt to provide general results that would guarantee \eqref{LMT} under different curvature bounds, such as $\ric\ge  b$, some $b \in \R$.
 

The proof of Theorem 1 is fairly short, based on the idea, used in [BM], that the lack of control of the size of the support of $u$ can be compensated by suitable truncation, an operation that preserves the class $W_0^{1,n}(\Omega)$. In [BM] the key tool was the P\'olya-Szeg\H{o} inequality $\|\nabla u^\#\|_n\le \|\nabla u\|_n$, where $u^\#$ is the Schwarz symmetrization of $u$, which is not always available on general manifolds. Instead, we make use of truncation inside a suitable level set of measure $\tau$. The point we wish to emphasize  here is the step from local to global, with a streamlined argument that works essentially in any setting where \eqref{LMT} and \eqref{Pn} make sense. As an example, on the Heisenberg group $\He^n$ one has a local Moser-Trudinger inequality as in \eqref{LMT}, for the horizontal gradient $\nabla_{\He^n}$, and  with a different, explicit, sharp constant $\alpha_Q$; this is due to Cohn-Lu [CL]. The argument presented in the proof of Theorem~1 works \emph{verbatim} in this sub-Riemannian  setting, thereby yielding the following result:

\smallskip
\begin{theorem} \label {3}   On the Heisenberg group $\He^n$,  given any open $\Omega\subseteq \He^n$ the sharp Moser-Trudinger inequality
\be\int_\Omega \exp_{ n-2}\big[\alpha_Q|u|^\nn\big] d\mu\le \Co,\qquad  u\in W_0^{1,n}(\Omega),\quad \|\nabla_{\He^n} u\|_n\le 1,\label{HMT}\tag{HMT}\ee
holds if and only if the Poincar\'e inequality $\|u\|_n\le \Co\|\nabla_{\He^n}u\|_n$ holds on $W_0^{1,n}(\Omega)$.
\end{theorem}

The validity of \eqref{LMT} on $\He^n$, which is basically \eqref{HMT} for  $\Omega$ bounded, has been proved by Cohn and Lu [CL], where the explicit value of $\alpha_Q$ has been computed.

\smallskip

We now discuss the validity of the global  Moser-Trudinger inequality under a stronger restriction on the norm of $u$. Given a complete, noncompact Riemannian manifold $M$, we say that the \emph{Ruf-Moser-Trudinger inequality holds on $W_0^{1,n}(M)$}, if for any  $\kappa>0$ there is  $C>0$ such that 
\be\int_M \exp_{ n-2}\big[\gamma_n|u|^\nn\big] d\mu\le C,\qquad  u\in W_0^{1,n}(M),\quad \|\nabla u\|_n^n+\kappa\|u\|_n^n\le 1.\label{RMT}\tag{RMT}\ee

By Observation 1, it's easy to check that  if \eqref{RMT} holds then the exponential constant must be sharp. Such  inequality was first proved in $\R^2$ by Ruf [R], and later extended to any $\R^n$ by Li-Ruf [LR]. The result was later extended to $\He^n$ in [CL], and to higher order Sobolev spaces in [FM1], [MasSani], [MQ].
As it turns out the validity of \eqref{LMT} is sufficient in order to have \eqref{RMT}:

\smallskip
 \begin{theorem}\label{4}  Let  $(M,g)$ be an  $n$-dimensional complete  Riemannian manifold. If  the local Moser-Trudinger inequality  \eqref{LMT} holds on $W_0^{1,n}(M)$,  then the Ruf-Moser-Trudinger inequality \eqref{RMT} holds on $W_0^{1,n}(M).$ 

\end{theorem}

The proof of Theorem \ref{4} is a very small modification of the proof of Theorem \ref{1}, and it can be traced back to the original paper by Ruf [R], who used  P\'olya-Szeg\H{o} inequality and truncation over suitable balls, or to Lam-Lu [LL1], where the authors use a similar truncation over  suitable level sets. The same argument has also been used on the Heisenberg group in  [LL2]. In [YSK] the authors derived \eqref{RMT} on Hadamard manifolds in a similar way, but, as we mentioned earlier, there was a gap in their proof of \eqref{LMT}. As an immediate  consequence of theorems \ref{2} and \ref{4} we obtain a version of their result which is valid under a  stronger curvature assumption:

\begin{cor}\label{5}  The sharp Ruf-Moser-Trudinger inequality \eqref{RMT} holds on any  $n$-dimensional complete, simply connected,   Riemannian manifold $(M,g)$ with sectional curvature $K$  satisfying $-b^2\le K\le 0$. 
\end{cor}

In [K], Krist\'aly proved that (RMT) holds on Hadamard manifolds satisfying the Cartan-Hadamard conjecture, where, essentially,  one can use the P\'olya-Szeg\H{o} inequality just as in $\R^n$. Such conjecture has been proved in dimensions $n=2,3,4$,  and therefore on such manifolds (RMT) holds without the restriction $K\ge -b^2.$

Once again, the main point we wish to emphasize in this paper  is that the step from local to global, on $W_0^{1,n}$, is essentially always true, wherever such inequalities make sense, and it is based on truncation on suitable level sets, combined with some elementary algebraic inequalities. This idea is not applicable  to sharp inequalities for higher order Sobolev spaces $W_0^{\alpha,n/\alpha}$, since truncation does not preserve such spaces, if $\alpha>1.$

\medskip
Finally, we note that while on compact manifolds we are not able to establish \eqref{LMT} for small $\tau$, 
the sharp RMT is still  valid, even for higher order gradients, without imposing restrictions on the averages of $u$ (see [LiY], [do\'OY]).

\medskip
\ni{\bf Acknowledgemt.} The authors would like to thank Anne Gallagher  for helpful discussions and comments.
\medskip
\section{Proofs of Theorems 1,2,4}

\ni{\underbar{\emph {Proof of Theorem \ref{1}.}}

\medskip

 First, it's easy to check that for any $u\in W_0^{1,n}(\Omega)$
\be \frac{\gn^{n-1}}{(n-1)!}\|u\|_n^n\le \int_\Omega\exp_{n-2}\big[\gn|u|^{\frac n{n-1}}\big]d\mu\le \int_{\{|u|\ge1\}}e^{\sgn |u|^{\frac n{n-1}}}d \mu+e^{\sgn}\|u\|_n^n\label{expreg}\ee
from which it's clear that if \eqref{CMT} holds then \eqref{Pn} holds as well.

Viceversa, suppose that \eqref{Pn} holds, so that from \eqref{expreg} it's enough to prove that 
\be \int_{\{|u|\ge1\}}e^{\sgn |u|^{\frac n{n-1}}}d\mu\le C,\qquad u\in W_0^{1,n}(\Omega),\quad \|\nabla u\|_n\le 1.\label{MT1}\ee
Note that under our hypothesis $\mu(\{|u|\ge1\})\le \|u\|_n^n\le C_\Omega^n$.

Take any $u\in W^{1,n}_0(\Omega)$, such that $\|\nabla u\|_n\le1$. Without loss of generality, we can assume $u\ge0$.

Note that 
\be \mu( \{x:u(x)>u^*(\tau)\})\le \tau\le \mu(\{x:u(x)\ge u^*(\tau)\})\ee
and that $u^*(\tau)>0$ if and only if $\mu(\supp u)>\tau$, where $\supp u$ denotes the essential support of $u$.

If $\mu(\supp u)\le \tau$ then \eqref{MT1} follows from \eqref{LMT}, so let us assume that $\mu(\supp u)>\tau$. In this case we can find a measurable set $E_\tau$ such that 
\be\mu(E_\tau)=\tau,\qquad \{x:u(x)>u^*(\tau)\}\subseteq E_\tau\subseteq\{x:u(x)\ge u^*(\tau)\}.\ee
For the existence of such a set, see  [MQ, (3.11)] and the comments thereafter.

Let 
\be v=(u-u^*(\tau))^+=\bc u-u^*(\tau) & {\text {if}}\; u\ge u^*(\tau)\cr 0 & {\text {if}}\; u<u^*(\tau).\cr\ec\ee
Then, $v\in W_0^{1,n}(\Omega)$, $\supp v\subseteq E_\tau$, and
\be \|\nabla v\|_n^n=\int_{E_\tau}|\nabla u|^nd\mu\le 1-\int_{E_\tau^c}|\nabla u|^n d\mu.\label {v1}\ee

Using the inequality $(a+b)^p\le a^p \e^{1-p}+b^p (1-\e)^{1-p}$ ($a,b\ge0,p\ge1, 0<\e<1$) we have
\be e^{\sgn u^\nn}\le  e^{\sgn v^\nn \e^{-\un} }\,e^{\sgn u^*(\tau)^\nn  (1-\e)^{-\un}}.\label{exp}
 \ee
From \eqref{Pn}
 we obtain 
 \be\ba 0<u^*(\tau)^n&=u^*(\tau)^n\;\frac{\mu(E_\tau)}{\tau}\le\frac1\tau\;\|u-v\|_n^n\cr& \le \frac{C_\Omega^n}\tau\;\|\nabla(u-v)\|_n^n =\frac{C_\Omega^n}\tau\int_{E_\tau^c}|\nabla u|^n\le \frac{C_\Omega^n}\tau\; (1-\|\nabla v\|_n^n)\label{u1}\ea\ee
 which implies $\|\nabla v\|_n<1$.
If $\|\nabla v\|_n>0$, let $\e=\|\nabla v\|_n^n$, in which case 
using \eqref{LMT} and \eqref{u1}
\be \int_{E_\tau}e^{\sgn v^\nn \e^{-\un} }\le C_0,\qquad e^{\sgn u^*(\tau)^\nn  (1-\e)^{-\un}}\le C\label{u2}\ee
which implies   \eqref{MT1}, since $\{u\ge1\}\subseteq E_\tau\cup \{1\le u\le u^*(\tau)\}\subseteq E_\tau\cup \{1\le u\le C_\Omega^n\}$, and, as we mentioned earlier,    $\mu(\{u\ge1\})\le~C_\Omega^n.$ If  $\|\nabla v\|_n=0$, then $u\le u^*(\tau)\le C_\Omega\tau^{-1/n}$ a.e., and \eqref{MT1} follows.

\rightline\qed

\bigskip
\ni{\underbar{\emph {Proof of Theorem \ref{4}}.} }

\medskip
The proof is identical to the one of Theorem \ref{1}, with the exception that assuming $ \|\nabla u\|_n^n+\kappa\|u\|_n^n\le 1$ we can  replace estimate \eqref{v1} with

\be  \|\nabla v\|_n^n=\int_{E_\tau}|\nabla u|^nd\mu\le 1-\kappa \|u\|_n^n-\int_{E_\tau^c}|\nabla u|^n d\mu.\label {v2}\ee
and \eqref{u1} with
 \be 0<u^*(\tau)^n=u^*(\tau)^n\;\frac{\mu(E_\tau)}\tau \le\frac1\tau\;\|u-v\|_n^n\le\frac1\tau\; \|u\|_n^n \le\frac1{\kappa\tau}\;( 1- \|\nabla v\|_n^n).\label{u3}\ee

\rightline\qed

\bigskip
\ni{\underbar{\emph {First proof of Theorem \ref{2}}.} }

\medskip

Let us assume $-b^2\le K\le 0$ and $n\ge 3.$ To prove \eqref{LMT}  it's enough to derive the rearrangement estimates in \eqref{kstar1}, \eqref{kstar2}, for the kernel $|\nabla_y G(x,y)|$, where $G(x,y)$ denotes the positive, minimal Green  function for $-\Delta$ on $M$. Those estimates, in turn,  follow from the bound
 \begin{numcases} {|\nabla_y G(x,y)|\le}\omega_{n-1}^{-1} d(x,y)^{1-n}\big(1+H_1 d(x,y)\big) & {\text {if} } $\; d(x,y)\le 1$ \label {GH1}\\ H_2 d(x,y)^{2-n} & {\text {if} } $\; d(x,y)> 1$\label {GH2},\end{numcases}
where $H_1,H_2$ are constants independent of $x,y,$ combined with the volume comparison estimate $\mu(B(x,r))\le V_{-b^2} (r)$, the volume of a ball of radius $r$ in the $n-$dimensional space form of constant curvature $-b^2$. Estimate \eqref{GH1} appears in  [BS, (3.9)],  whereas \eqref{GH2}  follows from the comparison theorem $G(x,y)\le  c_2 d(x,y)^{2-n}$ 
(where $c_2|x-y|^{2-n}$ is the classical Newtonian potential on $\R^n$) [BS, (3.3)], combined with the gradient estimate  $|\nabla_y G(x,y)|\le G(x,y)\big(C_5 d(x,y)^{-1}+(n-1)b\big)$ [Li, Thm. 6.1], [BS, (3.14)].  Since $|\nabla_y G(x,y))|\le H_2$ for $d(x,y)>1$ we get that if  $s>H_2$ then 
\be\mu(\{x:|\nabla_y G(x,y)|>s\})= \mu(\{x: d(x,y)\le 1,\; |\nabla_y G(x,y)|>s\})\ee
and  \eqref{GH1} together with the volume comparison $\mu(B(x,r))\le V_{-b^2}(r)\le  C_3 r^n$ for $r\le 1$, yield \eqref{kstar1}, \eqref{kstar2}, for some $t_0>0$.
\rightline \qed

\bigskip
\ni{\underbar{\emph {Second proof of Theorem \ref{2}}.} }

\medskip
 An alternative proof of Theorem 2, for any $n\ge2$,  can be given by using the same integral representation used in [YSK], namely \eqref{rep} with $k(x,y)=\omega_{n-1}^{-1}d(x,y)^{1-n} |g(y)|^{-1/2}$, where  $|g|=\det(g_{ij})$, and where $g=(g_{ij})$ is the metric tensor in geodesic polar coordinates centered at $x$ (recall that  $\inj(M)=+\infty,$ since $K\le 0$). Indeed, as proved in [YSK, Lemma 3.2] we have $(k_x)^*(t)\le \gamma_n^{-\frac {n-1}n} t^{-\frac{n-1}n}$ for $t>0$, for all $x$, hence \eqref{kstar1}  holds. The missing estimates in [YSK] was \eqref{kstar2}, which can be obtained with the additional assumption that $K\ge -b^2$. Indeed, since $|g(y)|\ge 1$ (due to the standard Bishop  comparison under $K\le 0$), then  $k(x,y)\le \omega_{n-1}^{-1} d(x,y)^{1-n}$, and the volume comparison $\mu(B(x,r))\le V_{-b^2}(r)$, arguing as in the previous proof, gives $(k^y)^*(t)\le C t^{-\frac{n-1}n}$ for $0<t\le t_0$, some $t_0>0$,  and all $y\in M$, that is \eqref{kstar2} also holds.

\rightline\qed

\medskip

\begin{rk}\label{curvature} Note that when $n\ge 3$ the global uniform estimate $|\nabla_y G(x,y)|\le C d(x,y)^{2-n}$, for $d(x,y)\ge 1$,  was derived using the gradient estimate and the volume comparison estimate,  and the latter are  possible due to the curvature  bound  $K\ge -b^2$ (which implies $\ric\ge -(n-1)b^2$). When $n=2$ this argument fails. In fact, we are not even able to show that one can always find a Green function $G(x,y)$ such that $|\nabla_y G(x,y)|\le C$, for $d(x,y)\ge 1$, assuming only  $-b^2\le K\le 0$.

\end{rk}

\medskip\section{Further results and comments}

\ni\emph{\underbar{The Poincar\'e inequality on $W_0^{1,n}(\Omega)$}}

\medskip

In view of Theorems 1,2,3 it is natural to ask for conditions that would guarantee that in a  given open set $\Omega$ the Poincar\'e inequality holds on $W_0^{1,n}(\Omega)$. On $\R^n$ this can be done using different notions of ``inradius". We  summarize here some of  what is known for completeness, hoping that it would motivate similar work  on more general manifolds (Riemannian or sub-Riemannian). 

For an open set $\Omega\subseteq \R^n$ the \emph{inradius} of $\Omega$ is defined as 
\be\rho(\Omega)=\sup\{R>0:\; \exists x\in\R^n  {\text { such that }} \;B(x,R)\cap\Omega^c=\emptyset\}\label{inr}\ee
the \emph{strict inradius} of $\Omega$  is defined as
 \be\rho_s(\Omega)=\sup\{R>0:\;\forall\e>0, \exists x\in\R^n {\text { such that }} \;B(x,R)\cap\Omega^c\not\supseteq B(z,\e), \forall z\in\R^n\}\label{sinr}\ee
 the \emph{measure-theoretic inradius} of $\Omega$ is defined as
 \be\rho_m(\Omega)=\sup\{R>0:\; \forall\e>0, \exists x\in\R^n {\text { such that }} \;|B(x,R)\cap\Omega^c|<\e\}\label{minr}\ee
 where $|E|$ denotes the Lebesgue measure of $E$, 
 and finally, the \emph{capacitary inradius} of $\Omega$
 is defined as \def\cp{\text{cap}}
  \be\rho_c(\Omega)=\sup\{R>0:\; \forall\e>0, \exists x\in\R^n {\text { such that }} \;\cp\big(B(x,R)\cap\Omega^c\big)<\e\}\label{cinr}\ee
 where $\cp(E)$ denotes either the Newtonian capacity  (if $n\ge3$) or the logarithmic capacity ($n=2$) of a Borel set $E$. 

The strict inradius was introduced by Souplet in [S], while  the capacitary inradius in the present form was defined by Gallagher-Lebl-Ramachandran  [GLR], and Gallagher [G1].    The  notion of capacitary inradius was introduced first  by Maz'ya-Shubin [MazShu] in a sligthly different form, namely for any $\alpha\in(0,1)$ they defined
\be\rho_{c,\a}(\Omega)=\sup\big\{R>0:\;  \exists x\in\R^n {\text { such that }} \;\cp\big(\ol{B(x,R)}\cap\Omega^c\big)\le\a\,\cp\big(\ol{ B(x,R)}\big)\big\}.\label{cinrb}\ee
 Likewise, based on the work of Lieb [L], Maz'ya-Shubin [MazShu] also defined the measure-theoretic inradius for a given $\alpha\in (0,1)$ as 
    \be\rho_{m,\a}(\Omega)=\sup\{R>0:\;  \exists x\in\R^n {\text { such that }} \;|B(x,R)\cap\Omega^c|\le\a |B(x,R)|\;\}.\label{minra}\ee
 It is easy  to check that 
\be \rho_m(\Omega)=\inf_{0<\alpha<1} \rho_{m,\a}(\Omega),\qquad \rho_c(\Omega)=\inf_{0<\alpha<1}\rho_{c,\a}(\Omega)\ee and
 \be \rho(\Omega)\le \rho_c(\Omega)\le \rho_m(\Omega)\le \rho_s(\Omega),\ee
 where the second inequality follows from the known bounds $\cp (E)\ge C_n |E|^{\frac {n-2}n}$ ($n\ge 3$) and $\cp(E)\ge \pi^{-1/2}|E|^{\frac12}$ ($n=2$).
 
 It is well known, see e.g. [S], that if $\rho(\Omega)=+\infty$ then the Poincar\'e inequality on $W_0^{1,p}(\Omega)$  i.e.
   \be\|u\|_p\le C\|\nabla u\|_p,\qquad \forall u\in W^{1,p}_0(\Omega)\label {PP}\ee
cannot hold for any $p\ge1.$ On the other hand in [S]  Souplet proves that if $\rho_s(\Omega)<\infty$ then \eqref{PP} holds for all $p\ge1$   (the special case $p=2$ is due to Agmon [Ag]).

  In [L], Lieb proves, among other things,  that if $\rho_{m,\a}(\Omega)<\infty$ for some $\alpha\in (0,1)$, or, equivalently, if  $\rho_m(\Omega)<\infty$,  then \eqref{PP} holds for all $p\ge1.$ 
 In [MazShu], Maz'ya-Shubin prove that the validity of \eqref{PP} for $p=2$, when $n\ge3$,  is equivalent to the condition $\rho_c(\Omega)<\infty$ (for a different proof see [G1]), and the same result is proved in [GLR],  when $n=2.$ In a very recent preprint [G2] Gallagher proves that \eqref{PP} is equivalent to the finiteness of the  $p$-capacitary inradius, for any $p>1$. The reader can consult the references just cited for specific examples of domains whose inradii of different types are or are  not finite.

\smallskip

 In a different direction,  in [FM2] the authors considered domains satisfying the condition
\be \supess_{x\in \Omega} \int_1^\infty \frac{|\Omega\cap B(x,r)|}{r^{n+1}}dr<\infty,\label{subcrit}\ee
 which were named ``Riesz subcritical'', in the sense that the rearrangement of the Riesz kernel $|x-y|^{\a-n}$ w.r. to $y\in \Omega$ is in $L^{\frac n{n-\a}}(1,\infty)$, for any $\a<n$, i.e. it decays better than $t^{-{\frac{n-\a}n}}$,  the rearrangement of $|x-y|^{\a-n}$ in  $ \R^n$. As proved in [FM2],  condition \eqref{subcrit}  implies - in fact \emph{it is equivalent to} - the following sharp  Adams inequality for the Riesz potential $I_\a f(x)=|x-y|^{\a-n}\star f$:
\be\int_\Omega \exp_{\llceil\frac n\a -2\rrceil }\bigg[\frac n{\omega_{n-1}}|I_\a f |^{\frac n{n-\a}}\bigg] dx\le C\big(1+\|I_\a f\|_{n/\a}^{n/\a}\big),\qquad  f\in L^{\frac n\a}(\Omega),\quad\|f\|_{n/\a}\le1,\label{rsub}\ee
which in turn implies - but it is \emph{not equivalent to} - the sharp Moser-Trudinger  inequality
\be\int_\Omega \exp_{\llceil\frac n\a -2\rrceil }\Big[\gamma_\a | u |^{\frac n{n-\a}}\Big] dx\le C\big(1+\|u\|_{n/\a}^{n/\a}\big),\qquad  u\in W_0^{\a,\frac n\a}(\Omega),\quad\|\nabla^\alpha u\|_{n/\a}\le1,\label{rsub1}\ee
where $\nabla^\alpha=(-\Delta)^{\frac\a2}$ for $\alpha$ even,  $\nabla^\alpha=\nabla (-\Delta)^{\frac{\a-1}2}$ for $\alpha $ odd ($\alpha<n$), and where  $\gamma_\alpha$ is the critical constant in the sharp Adams-Moser-Trudinger inequality on $\R^n$,  as in [A]. For the proof of these inequalities see [FM2], where it was  also observed that if in addition to \eqref{subcrit} the Poincar\'e inequality \eqref{Pn} holds, then  $\|u\|_{n/\a}\le C\|\nabla^\a u\|_{n/\a}$ also holds on $W_0^{\a,\frac n\a}(\Omega)$, and we can replace the RHS of \eqref{rsub1} by a constant (in particular \eqref{CMT} holds). As is turns out however \eqref{subcrit} actually implies that $\rho_m(\Omega)<\infty$, so \eqref{PP} holds for any $p\ge1.$ Thus, \eqref{subcrit}  gives another sufficient geometric condition for the validity of the Poincar\'e inequality. An example of $\Omega$ which satisfies \eqref{subcrit} but has $\rho_s(\Omega)=+\infty$ was given in [FM2, page 26].
\smallskip
\begin{prop} If $\Omega\subseteq \R^n$ is Riesz-subcritical, i.e. if \eqref{subcrit} holds, then $\rho_m(\Omega)<\infty$.\end{prop}

\ni{\bf Proof.} First, observe that if 
\be\rho_m'(\Omega)=\sup\{R>0:\; \forall\e>0, \exists x\in\Omega {\text { such that }} \;|B(x,R)\cap\Omega^c|<\e\}\label{minrp},\ee
then  clearly $\rho_m'(\Omega)\le\rho_m(\Omega)$, and  $\rho_m(\Omega)=+\infty$ if and only if $\rho_m'(\Omega)=+\infty$. Indeed, suppose that $\rho_m(\Omega)=+\infty$. Then, for each $R>1$ and $\e<|B(0,1)|:=C_n$ there is $x\in\R^n$ such that $|B(x,2R)\cap\Omega^c|<\e$. Hence 
$|B(x,2R)\cap\Omega|\ge C_n (2R)^n-\e$ and $|B(x,R)\cap\Omega|\ge C_n R^n-\e>0$. Therefore, there is $z\in B(x,R)\cap\Omega$, and $|B(z,R)\cap\Omega^c|\le |B(x,2R)\cap\Omega^c|<\e$, for all $R>1$ and  all $\e<C_n$, which is enough to imply that $\rho_m'(\Omega)=+\infty$.

Suppose then that $\rho_m'(\Omega)=+\infty$. Then, for each $R>1$ there is $x_R\in\Omega$ such that $|B(x_R,R)\cap \Omega^c|<1$. Hence for $1\le r\le R$, we have 
\be|\Omega\cap B(x_R,r)|=|B(x_R,r)|-|\Omega^c\cap B(x_R,r)|\ge |B(x_R,r)|-|\Omega^c\cap B(x_R,R)|\ge C_n r^n-1\ee
\be \supess_{x\in \Omega} \int_1^\infty \frac{|\Omega\cap B(x,r)|}{r^{n+1}}dr\ge \int_1^R \frac{C_n r^n-1}{r^{n+1}}dr\ge C_n\log R-\frac1n\to+\infty\ee
as $R\to+\infty$. 

\rightline \qed

\bigskip\ni\emph{\underbar{LMT and the volume of balls}}

\medskip
Below $B_x(r)$ denotes the geodesic ball with center $x$ and radius $r$.

\begin{prop}\label{lmtballs} If \eqref{LMT} holds then for each $\e>0$ there is $C_1=C(n,\e,C_0)>0$ such that
\be \mu(B(x,r))\ge\min\big\{\tau, C_1 r^{n+\e}\big\},\qquad\forall x\in M,\;\forall r>0.\label{balls}\ee
\end{prop}

The proof is similar to the one given by Y. Yang [Y, Lemma 3.1], which was itself a modification of an argument given by Carron [H, Lemma 3.2].  We note that in [Y, Prop 2.1], Yang showed that if \eqref{RMT} holds then $\inf_x \mu(B(x,r))\ge C>0$ for all $r$. Since \eqref{LMT} implies \eqref{RMT} we could just deduce the non collapsing  volume property from [Y, Prop 2.1]. However we give here a more direct proof, which also yields a slightly better lower bound than the one in [Y],  for small $r$.

\smallskip
\ni {\bf Proof.} Fix $x\in M$ and $r>0$. For simplicity let $B(r)=B_x(r)$, and  define
\be\psi(y)=\bc r-d(x,y) & \text{if} \;d(x,y)\le r \\ 0 & \text{if} \;d(x,y)>r. \ec \label{f}\ee
Then $\psi\in W_0^{1,n}(M)$ and $|\nabla \psi|=1$ a.e. on $B(r)$. Letting $u=\psi/\|\nabla \psi\|_n=\psi/\mu(B(r))^{1/n}$, if $\mu(B(r))<\tau$, then the  validity of \eqref{LMT} implies
\be\ba 1&\ge C_0^{-1}\int_{B(r)} \exp\big[\gamma_n u^{\frac n{n-1}}\big] d\mu\ge C_0^{-1}\int_{B(r/2)} \exp\big[\gamma_n u^{\frac n{n-1}}\big] d\mu\cr& \ge C_0^{-1} \exp\bigg[\gamma_n \frac{(r/2)^{\frac n{n-1}}}{\mu(B(r))^{\frac1{n-1}}}\bigg]\mu(B(r/2))\ge C_0^{-1} \frac{\gamma_n^q}{q!}  \frac{(r/2)^{\frac {nq}{n-1}}}{\mu(B(r))^{\frac q{n-1}}}\mu(B(r/2))
\ea\ee
Hence
\be \mu(B(r))\ge C(n,q) \Big(\frac r2\Big)^n \mu(B(r/2))^{\frac{n-1}q}\ee
where $C(n,q)=\Big(C_0^{-1} \frac{\gamma_n^q}{q!}\Big)^{\frac{n-1}q}$. By iteration it's easy to see that 
\be  \mu(B(r))\ge \big(C(n,q) r^n\big)^{\sum_0^m\left(\frac{n-1}q\right)^k}
\Big(\frac1{2^n}\Big)^{\sum_1^{m+1}k\left(\frac{n-1}q\right)^{k-1}}
 \mu(B(r/2^m))^{\left(\frac{n-1}q\right)^m}.
\ee
For any given $q>n-1$ we can let $m\to+\infty$, and using that $\mu(B(r))\sim \frac{\omega_{n-1}}n r^n$ we get $ \mu(B(r/2^m))^{\left(\frac{n-1}q\right)^m}\to1$ and 
\be \mu(B(r))\ge C(n,q)^{\frac q{q-n+1}} \Big(\frac1{2^n}\Big)^{\frac {q^2}{(q-n+1)^2}} r^{\frac{nq}{q-n+1}}, \ee
from which the result follows by taking $q$ large enough.

\rightline\qed

We note that in the above argument the specific value of $\gamma_n$ does not play any role. Hence the conclusion  of Proposition \ref{lmtballs} holds  regardless of the exponential constant in \eqref{LMT}.

\bigskip\ni\emph{\underbar{LMT on  $S^1\times\R$}}

\medskip
\def\p{\partial}

The flat 2-dimensional cylinder $M=S^1\times \R$ provides a simple example of a Riemannian manifold which is   complete, connected,  not simply connected, with  0 sectional curvature, and with  sub-Euclidean growth ($\mu(B(x,r))\sim C r$ for large $r$), such that \eqref{LMT} and \eqref{RMT} hold, but neither the Poincar\'e inequality nor \eqref{CMT} hold. 

To see this,  in the coordinate chart $\{(x,\theta): x\in\R,\; -\pi<\theta<\pi\}$ a Green function with pole $(0,0)$ for $\Delta$ is given by the explicit formula
\be G(\theta,x)=-\frac 1{4\pi}\log(\cosh x-\cos\theta)\ee
and clearly a Green function at pole $(\theta_0,x_0)$ is given by $G(\theta-\theta_0,x-x_0)$ (see [E, Ex. 2.5] for related discussions).

Clearly, if  $\nabla G=G_\theta\p_\theta+ G_x \p_x$ then 
\be(G_\theta,G_x)=-\frac1{4\pi} \Big(\frac{\sin \theta}{\cosh x-\cos\theta},\frac{\sinh x}{\cosh x-\cos\theta}\Big)\ee
and without much effort one sees that there is   $C_4>0$ such that 
\be\left||\nabla G(\theta,x)|-\frac 1{2\pi} (\theta^2+x^2)^{-1/2}\right|\le C_4,\qquad |x|\le 1,\; |\theta|<\pi\ee

 and
  \be\left||\nabla G(\theta,x)|^2-\frac1{16\pi^2}\right|\le \frac{1}{8\pi^2}\;\frac{1}{\cosh x-1},\qquad x\neq0,\; |\theta|<\pi\ee
   \be|\nabla G(\theta,x)|\le1,\qquad |x|\ge 1,\; |\theta|<\pi.\ee
   From the above estimates it is straightforward to check that \eqref{kstar1}, \eqref{kstar2} hold and therefore \eqref{LMT} and \eqref{RMT} hold on $S^1\times\R$. 
   
   The fact that both the Poincar\'e inequality and \eqref{CMT} do not hold on $S^1\times \R$ can be verified by a simple dilation argument in the $x$ variable, considering the family $u_r(\theta, x)=\phi(r x)$, where $\phi:\R\to\R$ is any nonzero, smooth, compactly supported function. Alternatively, one can use the family of functions defined in \eqref{f}.

\begin{table}[h]

\setlength{\tabcolsep}{24pt} 
\begin{tabular}{@{}p{0.5\linewidth}p{0.5\linewidth}@{}}
\textbf{Luigi Fontana} & \textbf{Carlo Morpurgo} \\
Dipartimento di Matematica ed Applicazioni & Department of Mathematics \\
Universit\'a di Milano-Bicocca & University of Missouri \\
Milan, 20125 & Columbia, Missouri 65211 \\
Italy & USA \\
\texttt{luigi.fontana@unimib.it} & \texttt{morpurgoc@umsystem.edu} \\\\
\textbf{Liuyu Qin} & \\
Department of Mathematics and Statistics & \\
Hunan University of Finance and Economics & \\
Changsha, Hunan & \\
China & \\
\texttt{Liuyu\_Qin@outlook.com} & \\
\end{tabular}
\end{table}
\end{document}